\documentclass[12pt]{amsart}
 
\usepackage{amsmath,amsthm,amssymb}
\usepackage{graphicx,psfrag,epsfig}

\begin{document}
\title{Diophantine properties of elements of ${\bf{SO(3)}}$.}
\author{V.~Kaloshin, I.~Rodnianski} 
\address {Department of Mathematics,
Princeton University, Princeton, NJ 08544}
\thanks{The first author is partially supported by 
the Sloan Dissertation Fellowship and the American Institute 
of Mathematics Five-year Fellowship}
\email{kaloshin@math.princeton.edu}
\email{irod@math.princeton.edu}
\maketitle
\markboth{Diophantine Rotations}{V. Kaloshin \& I. Rodnianski}       
 
\theoremstyle{plain}
\newtheorem{Thm}{Theorem}
\newtheorem{Def}[Thm]{Definition}
\newtheorem{Lm}{Lemma}
\newtheorem{Prop}{Proposition}
\newtheorem{Corollary}{Corollary}
\newtheorem{Rem}{Remark}

\def\bdef{\begin{Def}}
\def\endef{\end{Def}}
\def\bthm{\begin{Thm}}
\def\ethm{\end{Thm}} 
\def\blm{\begin{Lm}}
\def\elm{\end{Lm}}
\def\brm{\begin{Rem}}
\def\erm{\end{Rem}}
\def\beq{\begin{eqnarray}}
\def\eneq{\end{eqnarray}}
\def\Cal{\mathcal}
\def\th{\theta}
\def\Bbb{\mathbb}
\def\dt{\delta}
\def\eps{\epsilon}
\def\al{\alpha}
\def\bt{\beta}
\def\gm{\gamma}
\def\td{\tilde}
\def\R{\Bbb R}
\def\Z{\Bbb Z}
\def\T{\Bbb T}
\def\I{\Cal I}
\def\CR{\Cal R}
 
\section{Introduction}
The classical result of metric number theory on Diophantine 
properties of numbers says the following: for any $\eps>0$ and 
a.e. $\al \in \R$ the map $n \al (\mod 1)$ has a constant $C=C(\al)>0$ 
such that $n \al (\mod 1)>C|n|^{-1-\eps}$ for every integer $n$ {\cite{Kh}}. 

Diophantine properties of numbers arise in various
problems in metric number theory {\cite{Kh}}, smooth dynamical 
systems, holomorphic dynamics {\cite{HK}}, KAM theory {\cite{La}},
and others.

Generalizations of the metric number theory led to the development
of the theory of simultaneous Diophantine approximations and even
Diophantine approximations on manifolds. In the latter case consider
manifold $\,M\subset \Bbb R^n\,$ defined by $\,n\,$ analytic 
functions $\,f_1,..,f_n: U\subset \Bbb R^d\to \Bbb R$, 
$\,M=\{{\bf{f}}(x):\,x\in U\}$. Assume that functions
$\,1,f_1,..,f_n\,$ are linearly independent over $\,\Bbb R$.
One of the central questions of the theory is the following 
conjecture made by Sprind\u zuk in 1980 and recently proved 
by D. Kleinbock and G. Margulis {\cite{KM}}:

Any manifold $\,M\subset \Bbb R^n\,$ of the above type is extremal,
i.e., for almost all $\,{\bf{y}}\in M\,$ and any $\,\epsilon >0\,$ 
there exists a positive constant $D(y)$ such that for all 
$\,{\bf{q}}\in \Bbb Z^n\,$ and $\,p\in\Bbb Z\,$
\beq \label{dio}
|{\bf{q\cdot y}} + p|\ge \frac {D(y)}{\|{\bf{q}}\|^{n(1+\epsilon)}}.
\eneq
Here $\,{\bf{q\cdot y}}=\sum_{i=1}^n q_i y_i\,$ and
$\,\|{\bf{q}}\|=max_{1\le i\le n} |q_i|$.

In fact, Kleinbock-Margulis prove even a stronger statement that
$\,M\,$ is strongly extremal (approximation in the sense of
(\ref{dio}) is replaced by the notion of 
{\it{multiplicative approximation}}). The proof is based on 
the correspondence between approximation properties of number 
$\,{\bf{y}}\in\Bbb R^n\,$ and behavior of certain orbits in 
the space of unimodular lattices in $\,\Bbb R^{n+1}$.
 

The analogue of the Diophantune property can be also formulated in the
noncommutative setting. As far as we know very little is known in this
case. However, some intuition has already been developed for the group
$\,SU(2) (SO(3))$. We say that $\,g_1,..,g_k\in SU(2)\,$ are Diophantine
if there exists a positive constant $\,d(g_1,..,g_k)\,$ such that
for $\,n\ge 1\,$ and $\,W_n\,$ a word in $\,g_1,..,g_k\,$ of length $\,n\,$
$$
||W_n\,-\,Id||\,\ge\,d^{-n}.
$$
Our interest to the problem of Diophantine approximations on the group
$\,SO(3)\,$ stems mainly from the question formulated in the list of open
problems in the paper of A. Gamburd, D. Jakobson, and P. Sarnak (Problem 4):
{\it{The Haar generic elements 
$\,(g_1, g_2,..,g_k)\in SU(2)^k\,$ in the sense of measure are
Diophantine}} \cite{GJS}. The paper \cite{GJS} provides an elementary solution 
of Ruziewicz problem asserting that the Haar measure on $\,\Bbb S^2\,$ 
is the unique finitely additive $\,SO(3)\,$ invariant measure defined 
on Lebesgue sets.

In what follows it is more convenient for us to pass to the group $\,SO(3)\,$ 
and restrict our attention to the
case of two generators.
Consider a subgroup $\,F\,$ generated by two elements $\,A,B\in SO(3)$.
The group $SO(3)$ would have a Diophantine property if for almost
all rotations $\,A,B\in SO(3)\,$ in the sense of measure and all
reduced words $\,W_n\in F\,$ of length $\,n\,$ in 
$\,A,\ B,\ A^{-1},\ B^{-1}$, 
\beq \label{dioph} 
\| W_n\,-\,Id \|\,\ge\,D(A,B)^{-n}
\eneq
for some positive constant $\,D(A,B)$. The presence of the words of
the form $\,ABA^{-1}B^{-1}\,$ and like indicates that $\,F\,$ has to 
be a free subgroup. It is a classical fact that the set of elements 
$\,A,B\in SO(3)\,$ which do not generate a free subgroup is a countable
union of analytic sets of codimension one (see also Lemma \ref{nondeg} 
for an independent 
demonstration). To see this it is sufficient to establish the
existence of just one free subgroup of rank two. The first explicit 
construction of such a subgroup was given by Hausdorff in 1914 
in his work on Hausdorff-Banach-Tarski paradox. A free subgroup
$\,F\,$ of rank two in $\,SO(3)\,$ allows one to construct four
disjoint subsets of the sphere $\,\Bbb S^2\,$ such that after rotating 
these subsets by elements of $\,F\,$ one obtains two copies of 
$\,\Bbb S^2\,$ minus a countable set. Modulo the issue of the
countable set it follows that there is no finitely additive measure 
defined on {\it{all}} sets of $\,\Bbb S^2$. It also follows that any 
finitely additive $\,SO(3)\,$ invariant measure defined on Lebesgue 
sets is absolutely continuous with respect to the Lebesgue measure. 
The Ruziewicz problem is to show that any such measure necessarily 
coincides with the Lebesgue measure. In the general setting, the
problem is formulated for the finitely additive $\,SO(n+1)\,$
invariant measure on $\,\Bbb S^n$. It is interesting to note that 
in dimension one Banach provided a negative solution to the Ruziewicz 
problem. G. Margulis \cite{Ma} and D. Sullivan \cite{S} used Kazhdan
property (T) to give the positive answer in dimensions $\,n\ge 4$. 
For dimensions $\,n=2,3\,$ the affirmative solution had been given 
by V. Drinfeld \cite{D}.

The solution of Ruziewicz problem in dimensions $\,n\ge 2\,$ 
can be reduced to the problem of finding a free subgroup 
$\,F\in SO(n+1)\,$ with a {\it{spectral gap}} property \cite{R}.
Namely, consider the subspace 
$\,L^2_0(\Bbb S^n)=\{f\in L^2(\Bbb S^n):\,
\int\limits_{\Bbb S^n} f\,d\mu = 0\}$. 
Then $\,F\,$ is said to have a spectral gap property if there exists 
a positive constant $\,c\,$ such that for any 
$\,f\in L^2_0(\Bbb S^n)\,$ there exists an element $\,g\in F\,$ 
such that $\,\|f\circ g - f\|\ge c \|f\|$. After passing from 
$\,SO(3)\,$ to its double cover $\,SU(2)\,$ the above can
also be reformulated in terms of the spectra of the irreducible 
representations of $\,SU(2)\,$ restricted to the element 
$\,z=g_1+g_1^{-1}+..+g_k+g_k^{-1}$. Namely, let $\,\pi_N\,$ denote 
the irreducible representation of $\,SU(2)\,$ realized as 
a linear action on the space of homogeneous polynomials in two
variables of degree $\,N$. Define $\,\hat z(\pi_N)=
\pi_N(g_1)+\pi_N(g_1^{-1})+..+\pi_N(g_k)+\pi_N(g_k^{-1})\,$
to be an $\,(N+1)\times (N+1)\,$ matrix. Then we say that a subgroup 
$\,F\,$ generated by $\,g_1,..,g_k\,$ has a gap if   
$$
{\limsup}_{N\to\infty} \|\hat z(\pi_N)\| < \|z\|.
$$ 
A. Lubotzky, R. Phillips, and P. Sarnak construct explicit examples of 
elements $\,g_1,..,g_k\in SU(2)\,$ with $\,k\ge 3\,$ which generate 
a subgroup with a gap. For those
generators $\,\|\hat z(\pi_N)\|\le 2\sqrt{2k-1}<2k$ \cite{LPS1}.

Lubotzky-Phillips-Sarnak also show that the sequence of measures 
$\,\mu_N(z)\,$ associated with the eigenvalue distributions of 
$\,\hat z(\pi_N)\,$ has two accumulation points as $\,N\to \infty$. 
Namely, they prove that there exist two measures
$\,\nu^{\text{even}}(z)\,$ and $\,\nu^{\text{odd}}\,$
such that $\,\mu_{2N}(z)\to \nu^{\text{even}}(z)\,$ and 
$\,\mu_{2N+1}(z)\to \nu^{\text{odd}}(z)$. Moreover, the rate of the convergence 
depends on the Diophantine 
properties of the generators $\,g_1,..,g_k\,$ of $\,F$. In addition 
they show that a free subgroup generated by the elements 
$\,g_1,..,g_k\in SU(2)\,$ with algebraic entries is Diophantine.

In this paper we take a first step in an attempt to understand 
the Diophantine properties of the group $\,SO(3)$. We establish that 
almost all pairs of  rotations $\,(A,B)\in SO(3)\,$ generate subgroups that satisfy 
a {\it{weak}} Diophantine condition when the conjectured
exponent $\,n\,$ in (\ref{dioph}) is replaced by $\,n^2$.  
Although, the results below are stated for the rank two subgroups of
$\,SO(3)\,$ they can be easily generalized to include $\,SU(2)\,$ 
and higher number of generators. 
 
It follows from the pigeonhole principle and compactness of $SO(3)$
that an exponential estimate (not super-exponential) (\ref{dioph}) 
is the optimal one since the number of words of length $n$ grows 
exponentially with $n$. 
It is an easy exercise to show that for a Baire generic (residual) 
set of pairs $A,B \in SO(3)$ Diophantine condition is not satisfied. 
Therefore, the problem about Diophantine properties of elements 
of $SO(3)$ is another example of a property which fails on 
a Baire generic set, but holds on a set of full measure. 
Numerous examples
of this phenomena appear in dynamical systems and topology 
(see {\cite{O}}, {\cite{HSY}}, and {\cite{Ka}}).

As we mentioned above, in this paper we obtain the first result 
on Diophantine properties of elements of $SO(3)$. 
Consider $SO(3)$ with the Haar measure $\mu$ on it. 
We show that for an a.e. pair $(A,\ B) \in SO(3) \times SO(3)$ 
there is a constant $D>0$ such that for any $n$ and any 
word $W_n(A,B)$ of length $n$ in $A$ and $B$ we have
\beq \label{lowbound} 
\|W_n(A,B)\pm Id\|\geq D^{-n^2}.
\eneq

Let us describe the approach we use to prove the result and discuss the difficulties
which arise in the attempt to get an exponential lower bound as in (\ref{dioph}). 
Let $A, \ B \in SO(3)$ be two distinct elements, $k \in \Z_+$, and 
$W_n(A,B)$ be a word of length $n$ in $A$ and $B$. Denote by $\alpha$ 
and $\beta$ the angles of rotations of $A$ and $B$ respectively and 
by $\gamma$ the angle between the axes of $A$ and $B$. Without loss of generality we can 
assume that the axis of rotation of $A$, denote $v_A$, is the OX-axis 
in the ambient $\R^3$ and the axis of rotation of $B$, denote $v_B$, 
belongs to the $(x,y)$-plane forming angle $\gamma$ with $v_A$ in 
the clockwise direction. Notice any word $W_n(A,B)$ is uniquely
defined by a triple $(\alpha,\beta,\gamma)\in\T^3$. Denote 
$W_n(A,B)=W_n(\alpha,\beta,\gamma)$. Now consider the $3$-dimensional 
torus $\T^3$ as a parameter space with Lebesgue measure $m$.
It is clear that a set of full product Haar measure $\mu\times \mu$ on 
$SO(3)\times SO(3)$ corresponds to a set of full Lebesgue measure  
$m$ on $\T^3$.

The proof presented below is based on a standard Borel-Cantelli 
arguments. The rough sketch is as follows. Fix a word 
$W_n(\alpha,\beta,\gamma)$ of length $n$ in $A$ and $B$. The goal is 
to estimate the measure of the set of parameters $(\alpha,\beta,\gamma)\in\T^3$ 
for which $W_n(\alpha,\beta,\gamma)$ is at most
$D^{-n^2}$ away from $Id$. Let $\,m_n(D)\,$ be an upper bound for the measure
of the union of these sets over all words of length $n$. 
By Borel-Cantelli if $\sum_n m_n(D)< \infty$, then
for a.e. $(\alpha,\beta,\gamma)\in\T^3$ (\ref{lowbound})
holds for all except finitely many words. Increasing $D$
we satisfy those finitely many conditions and complete the proof.  

It turns out that a distance of $W_n(A,B)$ to  $Id$ can be 
represented as a trigonometric polynomial $P_n(\al,\bt,\gm)$ 
of degree $~n$ in $\al,\bt,$ and $\gm$ with integer coefficients.  
Fix $\bt=\bt^*$ and $\gm=\gm^*$ and consider measure of $\al$'s for 
which $P_n(\al,\bt^*,\gm^*)$ is $D^{-n^2}$-small. If 
a nontrivial $P_n(\al,\bt^*,\gm^*)$ with integer coefficients has a 
zero of order $n$ in $\al$ then measure 
$|\{\al: |P_n(\al,\bt^*,\gm^*)|<D^{-n^2}\}|$ can be as big as $~D^{-n}$. 
Suppose we can prove that $D^{-n}$ is an upper bound. Since, 
there are at most $4^n$ words $W_n(A,B)$ of length $n$ we obtain that 
the total ``bad'' measure of words of length $n$ is at most $(4/D)^n$ 
and is exponentially small for $D>4$.

One can think that the polynomial $P_n(\al,\bt^*,\gm^*)$
with a zero in $\al$ of high order corresponds to the fact
that the word $W_n(A,B)$ "sticks" in a neighborhood of $Id$ and leaves 
this neighborhood slowly as parameters $\al, \bt, \gm$ vary. This shows 
that a possible presence of high order degeneracies for the polynomial 
representing the
distance from a word $W_n(A,B)$ to $Id$ raises difficulties for
estimates of measure of a set where $W_n(A,B)$ is close to $Id$.
In particular, possible high degeneracies stand in the way of proving
the desired optimal result (\ref{dioph}).

In the last section we present a collection of words $W_n(A,B)$ 
of length $n$ for which polynomial $P_n(\al)$ does have a zero 
of order $\sqrt{n}$. This shows that it is indeed possible for a word $W_n(A,B)$ to 
"stick" in a neighborhood of $Id$. This degenerate collection is constructed 
using commutators $[A,B]=ABA^{-1}B^{-1}$. 
Degeneracies of high orders for trigonometric polynomials 
$P_n(\al,\bt,\gm)$ arising as a distance from a words
$W_n(\al,\bt,\gm)$ to $Id$ do occur.

\section{Statement of the result}

Let $A, \ B \in SO(3)$ be two distinct elements and 
$k \in \Z_+$. Denote $\I_m=(s_1,r_1,\dots,$\ $s_m,r_m)$
a set of $2m$ nonzero integers, $|\I_m|=\sum_p(|s_p|+|r_p|)$,
and $W_{\I_m}(A,B)=A^{s_1}B^{r_1} \dots A^{s_m} B^{r_m}$.
So, $W_{\I_m}(A,B)$ corresponds to the word defined by 
the multi-index $\I_m$. 
\bthm \label{main}
For any element $C \in SO(3)$ and $\mu\times \mu$-a.e.
pair $(A,B) \in SO(3)\times SO(3)$ there is a constant $D=D(A,B)>0$
such that
\beq
\min_{\{\I_m:\ |\I_m|=n\}} \|W_{\I_m}(A,B)-C\|\geq D^{-n^2}\ 
\ \ \textup{for all} \ \ \ n \in \Z_+.
\eneq

In other words, for $\mu$-generic choice of a pair $A$ and $B$,
all possible  words of length $n$ can not approximate
ahead given element $C$ better than $D^{-n^2}$.
The most interesting case when $C$ is the identity.
\ethm

Reformulate (\ref{main}) in a different form.
\bthm \label{main1}
For any element $C \in SO(3)$ and Lebesgue a.e.
$(\alpha, \beta, \gamma) \in \T^3$ there is a constant 
$D=D(\alpha, \beta, \gamma)>0$ such that
\beq \label{diophan}
\min_{\{\I_m:\ |\I_m|=n\}} 
\|W_{\I_m}(\alpha, \beta, \gamma)-C\|\geq D^{-n^2}\ 
\ \ \textup{for all} \ \ \ n \in \Z_+.
\eneq
\ethm
 
Fix a word $W_{\I_m}(\alpha, \beta, \gamma)$.
The idea of the proof is to show that outside of some
small measure set in $\T^3$ size of the derivative 
\beq
\|W_{\I_m}(\alpha, \beta, \gamma)'_{\alpha}\|^2=
D^\alpha_{\I_m}(\alpha, \beta, \gamma)
\eneq
is not too small. When the derivative with respect to $\al$ is not
too small the word $W_{\I_m}(\alpha, \beta, \gamma)$
varies  sufficiently fast with $\al$ and passes 
the ``dangerous'' $D^{-n^2}$-neighbor\-hood of 
the rotation $C$  sufficiently quickly. This implies
smallness of the ``prohibited'' set in the
parameter space  $(\alpha, \beta, \gamma)$.

Fix $n \in \Z_+$ and denote $\CR_n=\{\I_m:\ |\I_m|=n\}.$
Define 
\beq \label{one}
\begin{aligned}
\Phi_{\I_m}(D,C)=\{(\al,\bt,\gm)\in \T^3:\ &
\|W_{\I_m}(\alpha, \beta, \gamma)-C\|\leq D^{-n^2}\}\\
\Phi_n(D,C)=&\cup_{\I_m\in \R_n}\Phi_{\I_m}(D,C). 
\end{aligned}
\eneq

If for some $D^*>0$ we prove that
\beq \label{Borel}
\sum_{n=1}^\infty m\{\Phi_n(D^*,C)\}<\infty, 
\eneq
then for $m$-a.e. $(\alpha, \beta, \gamma)\in \T^3$
(resp. $\mu\times \mu$-a.e. $(A,B) \in SO(3)\times SO(3)$) 
there is $D=D(\alpha, \beta, \gamma) \geq D^*$
(resp. $D=D(A,B)$) such that (\ref{diophan})
is satisfied.

To estimate measure of $\Phi_n(D,C)$ we need to 
estimate measure of $\Phi_{\I_m}(D,C)$ for each
word $\I_m$ of length $n$, i.e. $|I_m|=n$.
Define the set of parameters, where the derivative
with respect to $\al$ is small
\beq \label{deriv}
\Phi^\al_{\I_m}=\{(\alpha, \beta, \gamma)\in \T^3:\ 
D_{\I_m}(\alpha, \beta, \gamma) \leq D^{-n^2/3}\}. 
\eneq

Denote $H(\R)$ the ring of quaternions 
$q=x_0 + i x_1 + j x_2 + k x_3, \ x_p \in \R$.
Let $\bar q=x_0 - (i x_1 + j x_2 + k x_3)$ and 
$N(q)=q \bar q.$ Denote $SH(\R)=\{q\in H(\R): N(q)=1\}$.
It is well-known that there is a representation of $SO(3)$ as $SH(\R)$
in the following form:
\beq \label{quat}
q=\cos \al + \sin \al (i v_1 + j v_2 + k v_3),
\eneq
where $\al$ is the angle of rotation
and a unit vector $(v_1,v_2,v_3)\in \R^3$
corresponds to an axis of rotation in the ambient
$\R^3$ of an element from $SO(3)$. 

\blm \label{odin} With the above notations
\beq
\|W_{\I_m}(\alpha, \beta, \gamma)''_{\al\al} \|^2\leq |\I_m|^4.
\eneq
\elm
{\it Proof}\ \  This follows from the quaternion representation
(\ref{quat}). Indeed, our choice of the ambient coordinate system gives
\beq \label{quat1}
\begin{aligned}
W_{\I_m}(\alpha, \beta, \gamma)=(\cos s_1\al + i \sin s_1 \al)
(\cos r_1\bt + \sin r_1 \bt(i \cos \gm + j \sin \gm)) \\
\dots (\cos s_m\al + i \sin s_m \al)
(\cos r_m\bt + \sin r_m \bt(i \cos \gm + j \sin \gm)).
\end{aligned}
\eneq
Differentiating this expression twice with respect to $\al$ gives
\beq
\|W_{\I_m}(\alpha, \beta, \gamma)''_{\al\al} \|^2 \leq
(\sum_{s=1}^k |s_p|)^4\leq |\I_m|^4. 
\eneq

\blm \label{nondeg}
The map $W_{\I_m}: SO(3)\times SO(3) \to SO(3)$
for a nontrivial word $\I_m$ is open. This, in particular,
implies that a pair of random elements of $SO(3)$ form 
a free group.
\elm

\brm
The conclusion of Lemma 2 is a well-known fact.
In particular, the statement that almost all subgroups in $\,SO(3)\,$ 
are free can be reduced to simply showing
that there exists {\it{a}} free subgroup in $SO(3)$. The latter is
a classical question which was solved positively first by F. Hausdorff
in 1914 \cite{Ha}. We present here a very explicit (constructive) 
independent proof of Lemma 2. 
\erm

{\it Proof}\ \   Consider representation (\ref{quat1}). 
To show that a trigonometric function is 
nontrivial with respect to, say $\al$, it is sufficient to establish that the 
highest frequency in $\al$ has a nonzero functional coefficient.
We shall compute this  functional coefficient,
namely, the coefficient in front the monomial 
$\exp (i \ \textup{sign}(s_m) \sum_{p=1}^m |s_p|\al)$. 
Notice that 
\beq \nonumber
\begin{aligned}\label{quat+}
&\textup{if}\ s>0\ e^{i s \al} (\cos r\bt + i \sin r\bt \cos \gm)=
(\cos r\bt + i &\sin r \bt \cos \gm) e^{i s \al}\\ 
&\label{quat-}
\textup{if}\ s<0\ e^{i s \al} j \sin r \bt \sin \gm =
j \sin r \bt \sin \gm e^{-i s \al}&
\end{aligned}
\eneq
Now we describe the procedure of permuting terms with 
$\al$ to the right and {\it particular terms  with $\bt$ 
and $\gm$} to the left so that after such  permutations 
the only term which has $\al$ is on the right end of the word  and equals
$\exp \left(i\ \textup{sign}(s_m) \sum_{p=1}^m |s_p|\al\right)$.

The first step of permutation: Consider the signs of $s_1$ 
and $s_2$. If they are different, then we change the sing 
of the $s_1$-term by choosing permutation (\ref{quat-}), 
otherwise, we choose (\ref{quat+}) in both cases
with $s=s_1$ and $r=r_1$. After the permutation the first
term with $\,\al$ from the left is 
$\exp \left(i\ \textup{sign}(s_2) \sum_{p=1}^2 |s_p|\al
\right)$.

The second step of permutation: Consider the signs of $s_2$ 
and $s_3$. Use the recipe of the first step.
The permutation gives the third term 
$\exp \left(i\ \textup{sign}(s_3) \sum_{p=1}^3 |s_p|\al\right)$ and so on.
Therefore, the only term which has 
$\exp \left(i\ \textup{sign}(s_m) \sum_{p=1}^m |s_p|\al\right)$ equals
\beq \nonumber
\begin{aligned}
 \prod_{\{p:\ s_p s_{p-1}>0\}}& (\cos r_p \bt + i\sin r_p \bt \cos \gm)\times
\\
 \prod_{\{p:\ s_p s_{p-1}<0\}}&  j \sin r_p\bt \sin \gm \quad
\exp \left(i\ \textup{sign}(s_m)\sum_{p=1}^m |s_p|\al\right).
\end{aligned}
\eneq 
This completes the proof.

\blm \label{dva} Let $|\I_m|=n$. Then
\beq
m\{\Phi_{\I_m}(D,C)\}\leq m\{\Phi_{\I_m}^\al(D,C)\}+4 D^{-n^2/3}n^4.
\eneq
\elm
{\it Proof}\ \ \  In the complement to the set 
$\Phi_{\I_m}^\al(D,C)$ we have estimates 
\beq \label{firstder}
\|W_{\I_m}(\alpha, \beta, \gamma)'_{\al} \|^2 \geq D^{-n^2/3}
\quad 
\|W_{\I_m}(\alpha, \beta, \gamma)''_{\al \al} \|^2 \leq n^4.
\eneq
For each pair $(\bt,\gm)\in \T^2_{\bt,\gm}$ split the circle 
$\T^1_\al$ into $\frac{D^{n^2/3}}{2n^4}$ intervals of equal 
length. Choose one interval and denote it by $I$. If there is a point 
in $(\alpha^*, \beta, \gamma) \in I$ which belongs to the complement
of $\Phi_{\I_m}^\al(D,C)$, then by the Taylor formula along with 
(\ref{firstder}) for each point in $I$ we have
 \beq
\|W_{\I_m}(\alpha, \beta, \gamma)'_{\al} \|^2 \geq \frac{D^{-n^2/3}}{2}.
\eneq
Therefore, the Taylor formula implies that measure of 
$\al \in I$ such that 
\beq
\|W_{\I_m}(\alpha, \beta, \gamma)-C\|\leq D^{-n^2}
\eneq
is at most $2 D^{-2n^2/3}$.  
Collecting all segments and applying Fubini's theorem
we complete the proof. 

Denote 
\beq
\Phi^\al_n(D)=\cup_{\I_m \in \CR_n} \Phi_{\I_m}^\al(D).
\eneq
Lemma \ref{dva} reduces a proof of (\ref{Borel}) to a proof
of 
\beq \label{Borel1}
\sum_{n=1}^\infty m\{\Phi^\al_n(D^*)\}<\infty. 
\eneq
We prove the convergence next.

\blm \label{tri} For any word $\I_m$ of length $n$ ($|\I_m|=n$)
there is a  polynomial
$P_{\I_m}(x_\al, y_\al,$\ $x_\bt, y_\bt, x_\gm, y_\gm)$
of degree $2n+2m$ with integer coefficients such that 
\beq
\quad \|W_{\I_m}(\alpha, \beta, \gamma)'_{\al} \|^2 =
P_{\I_m}(\cos \al, \sin \al, \cos \bt, \sin \bt, \cos \gm, \sin \gm).
\eneq
\elm

{\it Proof}\ \ Consider the quaternion representation 
(\ref{quat1}) differentiate it and take the sum of squares
of components. Then express $\cos s_p \al$ and $\sin s_p \al$ 
(resp. $\cos r_p \bt$ and $\sin r_p \bt$) as polynomials in
 $\cos \al$  and $\sin \al$ (resp. $\cos \bt$  and $\sin \bt$).
This gives a polynomial 
$P_{\I_m}(\cos \al, \sin \al, \cos \bt, \sin \bt, \cos \gm,
\sin \gm)$ with integer coefficients since all operations are with 
integer-coefficient trigonometric expressions. 
 
The main idea is that {\it a polynomial with integer coefficients
can not be small on a set of large measure}. In our notations
for $|\I_m|=n$ 
\beq \nonumber
\Phi^\al_{\I_m}(D)=\{(\al, \bt, \gm) \in \T^3:\ 
P_{\I_m}(\cos \al, \sin \al, \cos \bt, \sin \bt, \cos \gm,
\sin \gm)\leq D^{-n^2/3}\}.
\eneq

The following result for polynomials in one variable proved in the paper 
of S. Dani and G. Margulis 
\cite{DM}. For more general statements in this 
direction see also Kleinbock-Margulis \cite{KM}.
          
\blm \cite{DM, KM}
Let $\,F(x)\,$ be a polynomial of degree $\,\le n$. 
Denote $\|F\|_B:=\max_{x\in B} |F(x)|$. Then
for any open interval $\,B\,$
\beq \nonumber
 m_1\{x\in B:\,|F(x)|\,\le\,\epsilon\,\}\,\le\, 2n(n+1)^{\frac 1n} \left (\frac 
\epsilon{||F||_B}\right )^{\frac 1n} m_1\{B\}.
\eneq
\elm

\section{Elimination of Variables and Reduction to the 1-dimensional Case}
There are several technical difficulties that complicate matters in 
our setup. We need to show that a certain polynomial in several
variables does not spend too much time in the neighborhood of zero. 
In addition, we have a trigonometric polynomial which means that some 
of the variables are dependent. To resolve the latter we apply 
the procedure known as elimination of variables described in 
Lemma \ref{A} of next section. The former problem is treated
with the multiple application of Lemma \cite{DM} each time reducing 
the number of variables.
 
The polynomial in question is 
$\,P_{\I_m}(\cos\alpha,\sin\alpha,\cos\beta,\sin\beta,\cos\gamma,\sin\gamma)$.
We need an estimate on the size of the set $\Phi^\alpha_{\I_m}(D)$, 
defined above. The above set has essentially the same measure as the set
\beq
\begin{aligned}
\Bbb K:=&\{(x_\alpha,x_\beta,x_\gamma)\in [-1,1]^3:\,\,\
P_{\Cal I_m}(x_\alpha,y_\alpha,x_\beta,y_\beta,x_\gamma,y_\gamma)-
\epsilon\,=\,0,\\
y_\alpha^2  +x_\alpha^2-&1 =0,\,\,y_\beta^2+x_\beta^2-1=0,\,\,
y_\gamma^2+x_\gamma^2-1=0,
\,{\textup{for some}}\,\,\epsilon\,\leq\,D^{-n^2/3}\}.
\end{aligned}
\eneq
We will apply elimination of variables and Lemma \cite{DM} three times
in a row. First list properties of the polynomial  
$\,P_{\Cal I_m}(x_\alpha,y_\alpha,x_\beta,y_\beta,x_\gamma,y_\gamma)$.
\newline
$\,\bullet\,{\text{deg}}_{x_\alpha,y_\alpha,x_\beta,y_\beta,x_\gamma,y_\gamma} 
P\,\le\,2n$.
\newline
$\,\bullet\,P_{\Cal I_m}(x_\alpha,y_\alpha,x_\beta,y_\beta,x_\gamma,y_\gamma)=
\sum\limits_{l=0}^n p_l(x_\alpha,x_\beta,y_\beta,x_\gamma,y_\gamma) 
y_\alpha^l$,
\newline 
$\, \,\,|p_l|\leq H:=(2^n\,n)^2,\quad \forall l=0,..,n$.

Apply Lemma \ref{A} for the polynomials 
$\,\sum\limits_{l=0}^n p_l(x_\alpha,x_\beta,y_\beta,x_\gamma,y_\gamma) y_\alpha^l\,$
and $\,y_\alpha^2+x_\alpha^2-1\,$ with $\,s=r=2n\,$ and $\,H=(2^n n)^2$.
From the properties of the resultant 
$\,R_\epsilon(x_\alpha,x_\beta,y_\beta,x_\gamma,y_\gamma)\,$
defined in Lemma \ref{A} it follows that 
\beq
\begin{aligned}
\Bbb K\subset \{&(x_\alpha,x_\beta,x_\gamma)
\in [-1,1]^3:\,\,R_\epsilon(x_\alpha,x_\beta,y_\beta,x_\gamma,y_\gamma)=0,\\& 
\,y_\beta^2+x_\beta^2-1=0,\,\,y_\gamma^2+x_\gamma^2-1=0,
\,{\text{for some}} \,\,\epsilon\,\le\,D^{-n^2/3}\}.
\end{aligned}
\eneq
Using estimates (\ref{28}) we conclude that
\beq \label{20}
\begin{aligned}
 \Bbb K&\subset \{(x_\alpha,x_\beta,x_\gamma)
\in [-1,1]^3:\,\,R(x_\alpha,x_\beta,y_\beta,x_\gamma,y_\gamma)\,
\leq\,\delta,\\ 
y_\alpha^2 & +x_\alpha^2-1=0,\,y_\gamma^2+x_\gamma^2-1=0\}\\
\delta &:= D^{-n^2/3}(2^{2n}(2n\,H)\,+\,2^{2n}(2n\,H)^2 D^{-n^2/3}).
\end{aligned}
\eneq
Observe that $\,\delta\,$ is of the size $\,D^{-n^2}$.
Fix $\,(x_\beta,y_\beta,x_\gamma,y_\gamma)\,$ satisfying 
$\,y_\beta^2+x_\beta^2-1=0,\, y_\gamma^2+x_\gamma^2-1=0\,$ 
and apply Lemma \cite{DM} to the polynomial
$\,R(x_\alpha,x_\beta,y_\beta,x_\gamma,y_\gamma)\,$ with respect to $\,\al$. 
Let 
\beq
\Bbb K_{\beta,\gamma}:=\{x_\alpha\in [-1,1]:\,\,
(x_\alpha,x_\beta,x_\gamma)\in \Bbb K\}.
\eneq
It follows that
\beq \label{21}
m_1\{\Bbb K_{\beta,\gamma} \}\,\le\,16n(8n+1)^{\frac 1{8n}} \left (\frac 
{\delta}{||R(\cdot,x_\beta,y_\beta,x_\gamma,y_\gamma)||}\right )^{\frac 1{8n}}.
\eneq
Note that $\,m_1\,$ and $\,m_2\,$ denote one and two-dimensional Lebesgue
measures correspondingly.\newline
Define 
\beq
\begin{aligned}
\Bbb K^1:=\{&(x_\beta,x_\gamma)\in [-1,1]^2:\,\,
||R(\cdot,x_\beta,y_\beta,x_\gamma,y_\gamma)||\,\le\,
\delta^{\frac 12},\\
&y_\beta^2+x_\beta^2-1=0,\,\, y_\gamma^2+x_\gamma^2-1=0\}.
\end{aligned}
\eneq
The Fubini Theorem implies that 
\beq \label{22}
m\{\Bbb K\}\le\, 2 m_2\{\Bbb K^1\}+ 
m\Bigl\{\bigcup_{(x_\beta,x_\gamma)\not\in \Bbb K^1}\Bbb K_{\beta,\gamma}\Bigr\}.
\eneq
Observe also that by the Fubini Theorem and (21) the set 
$\,\bigcup_{(x_\beta,x_\gamma)\not\in \Bbb K^1}\Bbb K_{x_\beta,x_\gamma}\,$
obeys the following estimate on its size:
\beq \label{23}
m\Bigl\{\bigcup_{(x_\beta,x_\gamma)\not\in \Bbb K^1}\Bbb K_{x_\beta,x_\gamma}\Bigr\}\,
\le\,64n(8n+1)^{\frac 1{8n}}  
\delta^{\frac 1{16n}}.
\eneq
 
To estimate the size of the set $\,\Bbb K^1\,$ we employ the conclusions
of the second part of Lemma \ref{A}. 
Define $\,P^1_{\Cal I_m}(x_\beta,y_\beta,x_\gamma,y_\gamma)\,$ from 
the resultant $\,R(x_\alpha,x_\beta,y_\beta,x_\gamma,y_\gamma)\,$ as in (\ref{29}):
\beq
P^1_{\Cal I_m}(x_\beta,y_\beta,x_\gamma,y_\gamma):=\,(16n)!\int\limits_{-1}^1
|R(x_\alpha,x_\beta,y_\beta,x_\gamma,y_\gamma)|^2\,dx_\alpha.
\eneq
The constant in front of the integral is introduced so that the resulting polynomial
is still a polynomial with integer coefficients.
Clearly,
\beq \label{24}
\begin{aligned} 
\Bbb K^1\subset\{&(x_\beta,x_\gamma)\in [-1,1]^2:\,
|P^1_{\Cal I_m}(x_\beta,y_\beta,x_\gamma,y_\gamma)|
\,\le\,2(16n)!\delta,\\
&y_\beta^2+x_\beta^2-1=0,\,\, y_\gamma^2+x_\gamma^2-1=0\}.
\end{aligned}
\eneq
Combining estimates (\ref{22}), (\ref{23}), and (\ref{24})  we
conclude that there exist positive constants $\,C_1,\rho_1\,$
such that
\beq
\begin{aligned}
m\{\Phi^\alpha_{\Cal I_m}(D)\}\,\le\,2\,m_2\bigl\{&(x_\beta,x_\gamma)\in [-1,1]^2:\,\,
|P^1_{\Cal I_m}(x_\beta,y_\beta,x_\gamma,y_\gamma)|\,\le\,D^{-\rho_1n^2},\\
&y_\beta^2+x_\beta^2-1=0,\,\, y_\gamma^2+x_\gamma^2-1=0\bigr\}\,+\,C_1^{-n}.
\end{aligned}
\eneq
The problem is now reduced to a similar two-dimensional question.
We are in position to apply another round of Lemma \ref{A}
and Lemma \cite{DM}. Reiterate the arguments above
for the polynomial $\,P^1_{\Cal I_m}(x_\beta,y_\beta,x_\gamma,y_\gamma)\,$ 
with properties as described
in Lemma \ref{A}:\newline
$\,\bullet\,\,{\text{deg}}_{x_\beta,y_\beta,x_\gamma,y_\gamma} 
P^1_{\Cal I_m}\,\le\,16n$.\newline
$\,\bullet\,P^1_{\Cal I_m}(x_\beta,y_\beta,x_\gamma,y_\gamma)\,=\,
\sum_{l=0}^{16n} p_{1l}(x_\beta,x_\gamma,y_\gamma) y_\beta^l$, and 
\newline $\,\bullet\,\,\max_{x_\beta,x_\gamma,y_\gamma\in [-1,1]^3} 
|p_{1l}(x_\beta,x_\gamma,y_\gamma)|\le H_1:=((16n)!)^3 4^{2n+1} (2nH)^4$.\newline
Note that by a crude estimate for any positive 
$\,\epsilon\,$ and all sufficiently large $\,n$, 
$\,H_1\le 2^{n^{1+\epsilon}}$.

Define the resultant $\,R^1(x_\beta,x_\gamma,y_\gamma)\,$ of the polynomials
$\,P^1_{\Cal I_m}(x_\beta,y_\beta,x_\gamma,y_\gamma)\,$ and 
$\,y_\beta^2+x_\beta^2-1$.
$\,$\newline
We obtain 
\beq \label{25}
\begin{aligned}
m_2\{(x_\beta,x_\gamma)\in [-1,1]^2:\,\,
|P^1_{\Cal I_m}(x_\beta,y_\beta,x_\gamma,y_\gamma)|\,\leq\,D^{-\rho_1n^2},\\
y_\beta^2+x_\beta^2-1=0,\,y_\gamma^2+x_\gamma^2-1=0\}\,\leq\,2\,m_1\{\Bbb K^2\}+ 
m_2\left\{\bigcup_{x_\gamma\not\in \Bbb K^2}\Bbb K^2_{x_\gamma}\right\},\\
\Bbb K^2:=\{x_\gamma\in [-1,1]:\,\,\|R^1(\cdot,x_\gamma,y_\gamma)\|\leq 
\delta_1^{\frac 12},\,\,  y_\gamma^2+x_\gamma^2-1=0\},\\
\delta_1:=D^{-\rho_1n^2}(2^{32n}16n\,H_1+2^{32n}(16n\,H_1)^2 D^{-\rho_1n^2}),\\
m_2\left\{\bigcup_{x_\gamma\not\in \Bbb K^2}\Bbb K^2_{x_\gamma}\right\}\,
\leq\,4(16n)(4(8n)+1)^{\frac 1{64n}}  
\delta_1^{\frac 1{64n}}.
\end{aligned}
\eneq
Observe that $\,\delta_1\,$ is still of the size $\,D^{-n^2}$. Therefore,
there exist positive constants $\,C_2,\rho_2\,$ such that
\beq
\begin{aligned}
m\{\Phi^\alpha_{\Cal I_m}(D)\}\,\leq\,2\,m_1\{x_\gamma\in [-1,1]:\,\,
|P^2_{\Cal I_m}(x_\gamma,y_\gamma)|\,\leq\,D^{-\rho_2n^2},\\
y_\gamma^2+x_\gamma^2-1=0\}\,+\,C_2^{-n}\,+C_1^{-n},
\end{aligned}
\eneq
where the polynomial $\,P^2_{\Cal I_m}(x_\gamma,y_\gamma)\,$ is formed from the
resultant $\,R^1(x_\beta,x_\gamma,y_\gamma)\,$ as in (\ref{29}).
Finally, eliminating $\,y_\gamma\,$ and applying Lemma \cite{DM} we can find
a positive constant $\,\delta_2\,$ of the size $\,D^{-n^2}\,$ 
such that
\beq \label{26}
m\{\Phi^\alpha_{\Cal I_m}(D)\}\,\leq\,
\Bigl ( \frac {\delta_2}{\|R^2(\cdot)\|} \Bigr )^{\frac 1{256 n}}
+\,C_2^{-n}\,+\,C_1^{-n}.
\eneq
The resultant $\,R^2(x_\gamma)\,$ is a polynomial with integer
coefficients of degree at most $\,128n$. Therefore, 
$\,(256n)!\int\limits_{-1}^1 |R_2(x_\gamma)|^2\,dx_\gamma\,$ is a 
non-negative integer. If it is positive, the desired estimate
immediately follows from (\ref{26}). So  we need to make sure that 
$R_2(x_\gamma)$ is not identically zero.
 
The polynomial $R_2(x_\gamma)$ was obtained via combination of elimination of variables 
(forming the resultant) and integration as in (\ref{29}).
Certainly, integration can not produce the identically zero
polynomial from a nonzero one. Therefore, we need to justify
the ``non-degeneracy'' of elimination. The basic property of the resultant
$R[P_1,P_2](x)$ of two polynomials $P_1(x,y), P_2(x,y)$, defined
below in (\ref{resultant}), is that $R[P_1,P_2](x_0)$ equals $0$ if and 
only if for some $y\in \Bbb C$ we have $P_1(x_0,y)=P_2(x_0,y)=0$
({\cite{Mu},p.34}). In our case one of polynomials, say $P_2$, is 
$x^2+y^2-1$. If $R[P_1,P_2](x)\equiv 0$, then $\,x=\cos\alpha, y=\sin\alpha$, and
$P_1(\cos \al, \sin \al)$ 
vanishes on the open set $\al \in U \subset \R$, which 
implies that it is identically zero. This is in contradiction
with non-degeneracy of $W_{\I_m}(\al,\bt,\gm)$ (see Lemma \ref{nondeg}).

\section{An Auxiliary Lemma}

Let $\,P_1(x,y)\,=\,\sum_{l=0}^r p_{1l}(x)\,y^l\,$ and   
$\,P_2(x,y)\,=\,\sum_{l=0}^s p_{2l}(x)\,y^l\,$ be two polynomials in $\,y\,$
of degree $\,r\,$ and $\,s\,$ correspondingly. Define the {\it{resultant}}
$$
R[P_1,P_2](x):=\,\,\det\,\Bbb A
$$ 
of $\,P_1\,$ and $\,P_2\,$ as the determinant of the following 
$\,(r+s)\times (r+s)\,$ matrix
\beq \label{resultant}
\Bbb A:=\left(\begin{array}{cccccccc}
p_{1r}(x)&\hdots&p_{10}(x)&0&\hdots&\hdots&\hdots&0\\
0&p_{1r}(x)&\hdots&p_{10}(x)&0&\hdots&\hdots&0\\
\vdots&\vdots&\vdots&\vdots&\vdots&\vdots&\vdots&\vdots\\
\vdots&\vdots&\vdots&\vdots&\vdots&\vdots&\vdots&\vdots\\
0&\hdots&\hdots&\hdots&0&p_{1r}(x)&\hdots&p_{10}(x)\\
p_{2s}(x)&\hdots&\hdots&p_{20}(x)&0&\hdots&\hdots&0\\
0&p_{2s}(x)&\hdots&\hdots&p_{20}(x)&0&\hdots&0\\
\vdots&\vdots&\vdots&\vdots&\vdots&\vdots&\vdots\\
\vdots&\vdots&\vdots&\vdots&\vdots&\vdots&\vdots\\
0&\hdots&\hdots&0&p_{2s}(x)&\hdots&\hdots&p_{20}(x)
\end{array} \right)
\eneq

We formulate an auxiliary lemma
\blm \label{A}
Let $\,P(x,y,u,v)\,=\,\sum_{l=0}^r p_{l}(x,u,v)\,y^l$.
Assume that the coefficients $\,p_l(x,u,v)\,$ are polynomials
of $\,(x,u,v)\,$ of degree $\,\leq s\,$ with respect to each variable:
$\,{\text{deg}}_{x,u,v} p_l\,\leq\, s$.

Assume also that for some constant $\,H\geq 1\,$ there holds the following 
estimates
\beq \label{27}
\max_{x,u,v\in [-1,1]} |p_l(x,u,v)|\,\leq\,H,\qquad\forall l=0,..r.
\eneq
Form a resultant $\,R_\epsilon (x,u,v)\,$ of the polynomials
$\,P(x,y,u,v)-\epsilon\,$ and $\,y^2+x^2-1$.
Then

$\,\bullet\,$ If for some $\,y\,$ the polynomials 
$\,P(x,y,u,v)-\epsilon=y^2+x^2-1=0$, the resultant $\,R_\epsilon (x,u,v)=0$.

$\,
\bullet\,R_\epsilon =R + \epsilon R_1 + \epsilon^2 R_2$,

where $\,R(x,u,f)\,$ is the resultant of $\,P(x,y,u,v)\,$ and $\,y^2+x^2-1$,
and
\beq \label{28}
\max_{x,u,v\in [-1,1]} |R_i(x,u,v)|\,\le\,2^{r}(rH)^{2-i},\qquad
i=0,..,2.
\eneq
Define the following polynomial of $(u,v)$:
\beq \label{29}
P_1(u,v):=\,(4(s+r))!\int\limits_{-1}^1 (R(x,u,v))^2\,dx.
\eneq
Then

$\,\bullet \,{\text{deg}}_{u,v} P_1\,\le\,4(s+r)$.

The polynomial $\,P_1(u,v)\,$ can be written as

$\bullet\,P_1(u,v)\,=\,\sum_{l=0}^{4(s+r)} p_{1l}(u)\,v^k$,

and
\beq \label{30}
\begin{aligned}
\max_{u\in [-1,1]} |p_{1l}(u)|\,\leq\,
\frac{((4(s+r)-l)!)^3}{l!} 
4^{r+1}(rH)^4,\quad
\forall l=0,..,4(s+r).
\end{aligned}
\eneq
\elm
{\it Proof}\ \ \ 
The $\,(r+2)\times(r+2)\,$ matrix corresponding to the resultant of the polynomials 
$\,P(x,y,u,v)-\epsilon\,$ and $\,y^2+x^2-1$ has the form
\beq \label{31}
\Bbb A:=\left(\begin{array}{cccccc} 
1&0&x^2-1&0&\hdots&0\\
\vdots&\vdots&\vdots&\vdots&\vdots&\vdots\\ 
\vdots&\vdots&\vdots&\vdots&\vdots&\vdots\\
\vdots&\vdots&\vdots&\vdots&\vdots&\vdots\\
0&\hdots&0&1&0&x^2-1\\
p_{r}(\cdot)&\hdots&\hdots&\hdots&p_{0}(\cdot)-\epsilon &0\\
0&p_{r}(\cdot)&\hdots&\hdots&\hdots&p_{0}(\cdot)-\epsilon\\
\end{array}\right)
\eneq
Any solution $\,y\,$ of the system $\,P(x,y,u,v)-\epsilon = y^2+x^2-1 = 0\,$
produces a nontrivial kernel containing the vector $\,(y^{r+1},y^r,..,1)\,$
of the matrix $\,\Bbb A$. Therefore, if for fixed
$\,(x,u,v)$ such a $\,y\,$ exists, the resultant $\,R_\epsilon(x,u,v)\,$ vanishes.

The estimate (\ref{30}) is the only nontrivial remaining statement of this lemma.
Its proof is based on the application of the Markov inequality:
$$
\max_{x\in [-1,1]} | F^\prime (x) |\,\leq\,n^2\,\max_{x\in [-1,1]}| F(x) |
$$
for any polynomial $\,F\,$ of degree $\,n$.
It easily follows from (\ref{28}) and (\ref{29}) that 
$$
\max_{u,v\in [-1,1]}| P_1(u,v) |\,\leq\, (4(s+r))!\,4^{r+1}(rH)^4.
$$
The coefficient $\,p_{1l}(u)\,$ can be found from the identity
$$
p_{1l}(u)\,=\,\frac 1{l!}\frac{d^l}{dv^l} P_1(u,0)
$$
Using Markov's inequality for the polynomial $\,P_1(u,v)\,$ of degree
$\,4(s+r)\,$ $\,l\,$ times we conclude that
$$\aligned
|p_{1l}(u)|\,\leq\,\frac{((4(s+r)-l)!)^3}{l!}4^{r+1}(rH)^4,\quad
\forall l=0,..,4(s+r).
\endaligned
$$

\section{Degenerate Words}
In this section for $n=4^m \in Z_+$ we construct $\sqrt{n}$ 
words $W_n(\al,\bt,\gm)$ such that if $P_n(\al,\bt,\gm)$ is 
the polynomial of distance of $W_n(\al,\bt,\gm)$ to $Id$, 
defined above, then it has a zero of order $\sqrt{n}$ with respect 
to $\al$ at any point of the form $(0,\bt,\gm)$.

Recall that $W_n(\al,\bt,\gm)=W_n(A,B)$ is a word in $A$ and
$B$, defined by the angle of rotation $\al$ of $A$, the angle of rotation 
$\bt$ of $B$ and the angle $\gm$ between the axis of rotations of $A$ 
and $B$ (see the introduction). Denote by $[A,B]=ABA^{-1}B^{-1}$ 
the commutator formed by $A$ 
and $B$. The idea of the construction is the following remark:
For a sufficiently small $\al$ the angle of rotation of the commutator $[A,B]$  
is of order at most $\al^2$.
At most $\al^2$ because, if axis of $A$ and $B$ are $\al$-close,
then  $[A,B]$ has an angle of rotation of order at most $\al^3$.
This follows directly from the quaternion representation
(\ref{quat}).

Consider two rotations $\Cal A, \Cal B \in SO(3)$. 
Define a map 
\beq
\phi:{\Cal A\choose \Cal B} \mapsto 
{{\Cal A \Cal B \Cal A^{-1} \Cal B^{-1}} \choose 
{\Cal B \Cal A \Cal B^{-1} \Cal A^{-1}} },
\eneq
which maps a pair of rotations into a pair of
commutator rotations. Define
\beq
\phi:{\Cal A_{k+1}\choose \Cal B_{k+1}} \mapsto 
{\Cal A_k \Cal B_k \Cal A^{-1}_k \Cal B^{-1}_k \choose 
\Cal B_k \Cal A_k \Cal B^{-1}_k \Cal A^{-1}_k},
\eneq
where $\Cal A_0=A^{\pm 1}$ and $\Cal B_0=B^{\pm 1}$.
Notice that $\Cal A_1$ and $\Cal B_1$ are rotations by 
an angle of order at most $\al^{2}$ provided that 
$\al$ is sufficiently small. $\Cal A_2$ and $\Cal B_2$ 
are rotations by an angle of order at most $\al^{4}$,
and $\Cal A_k$ and $\Cal B_k$ are  rotations 
by an angle of order at most $\al^{2^k}$. Since there 
is freedom in choosing powers of $\Cal A_k^{\pm}$ and 
$\Cal B_k^{\pm}$ in the definition of $\phi$ it is easy 
to see that this construction gives at least $2^k$ words 
of kind $\Cal A_k^{\pm}$ and $\Cal B_k^{\pm}$. Note that
$\Cal A_k^{\pm}$ and $\Cal B_k^{\pm}$ are words of length $4^k$. 

Let $\rho$ be the golden mean. It is not too difficult to see
that after choosing $A^{\pm 1}$ and 
$B^{\pm 1}$ in an appropriate way inside of the commutators 
one can construct  a word 
$W_n(A,B)$ with a zero of order $n^{(\rho+1)/2}=2^{(\rho+1)k}$ 
at the point $\,(\al,\bt,\gm)=(0,\bt,\gm)$.
 
All degenerations described here occur in a neighborhood
of zero. It is an interesting question whether there are
zeroes of high order far away from the identity element in
$\,SO(3)$.

\medskip

{\it Acknowledgments:}\ \ We would like to thank 
Peter Sarnak for drawing our attention to the problem, stimulating discussions,
and encouragement. We also thank Dmitriy Jakobson for valuable suggestions.


\medskip  

\begin{thebibliography}{GNED}
\def\bitem#1{\bibitem[#1]{#1}}

\bitem{DM} Dani S., Margulis G. Limit distributions of orbits of
unipotent flows and values of quadratic forms. 
I. M. Gelfand Seminar, pp.91--137, Adv. Soviet Math., {\bf{16}},
Part 1, Amer. Math. Soc., Providence, RI, 1993;

\bitem{D} Drinfeld, V. Finitely additive measures on 
$\,\Bbb S^2\,$ and $\,\Bbb S^3$, invariant with respect to 
rotations. Func. Anal. and its Appl.
{\bf{18}}, pp.245--246, 1984;

\bitem{HK} Hasselblatt, B. Katok, A. Introduction to the modern theory 
of dynamical systems. Encyclopedia of Math and its App, {\bf{54}}. 
Cambridge University Press, Cambridge, 1995;  

\bitem{Ha} Hausdorff, F. Grunzuge der Mengenlehre, Leipzig, 1914

\bitem{HSY} Hunt, B. Sauer, T. Yorke, J. Prevalence: 
a translation-invariant "almost every" on 
infinite-dimensional spaces.
Bull. Amer. Math. Soc. {\bf{27}}, no. 2, pp.217--238, 1992;
{\bf{28}}, no. 2, pp.306--307, 1993.

\bitem{GJS} Gamburd A., Jakobson, D., Sarnak, P.
Spectra of elements in the group ring of $SU(2)$. 
J. Eur. Math. Soc. (JEMS), {\bf{1}}, no. 1, 
pp.51-85, 1999;

\bitem{Ka} Kaloshin, V. Some prevalent properties of 
smooth dynamical systems. Proc. of Steklov Math. Inst., {\bf{213}}, 
pp. 123--151, 1997; 

\bitem{Kh} Khintchine, A. Continued fractions. Translated by P.Wynn. 
P. Noordhoff, Ltd., Groningen 1963

\bitem{KM} Kleinbock D., Margulis, G. Flows on homogeneous spaces and
Diophantine approximation on manifolds. Ann. of Math. {\bf{148}},
no. 1, pp.339--360, 1998; 

\bitem{La} Lazutkin, V. KAM theory and semiclassical approximations to 
eigenfunctions. Ergeb. Math. Grenzgeb.(3), 24. Springer-Verlag,
Berlin, 1993;

\bitem{Lu} Lubotzky, A. Discrete groups, expanding graphs
and invariant measures. Progress in Mathematics 125, Birkh\"auser,
Basel, 1994;

\bitem{LPS1} Lubotzky A., Phillips R.,  Sarnak P.,
 Hecke operators and distributing points on the
sphere. I. Comm. Pure Appl. Math., {\bf{39}}, no. S, suppl., 
pp.149-186, 1986;

\bitem{LPS2} Lubotzky A., Phillips R.,  Sarnak P.,
 Hecke operators and distributing points on the
sphere. II. Comm. Pure Appl. Math., {\bf{40}}, no.4, pp.401-420,
1987.

\bitem{Ma} Margulis, G. Some remarks on invarinat means.
Monatschefte fur Mathematik, {\bf{90}}, pp.233--235, 1980;

\bitem{Mu} Mumford, D. Algebraic Geometry I, Complex Projective 
Varieties. Springer-Verlag, New York, 1976;


\bitem{O} Oxtoby, J. Measure and category.  Grad Texts in Math, {\bf{2}},
Springer-Verlag, New York-Berlin, 1980.  

\bitem{R} Rosenblatt, J. Uniqueness of invarinat means for measure 
preserving transformations. Trans. AMS, {\bf{265}}, pp. 623-636, 1981;

\bitem{S} Sullivan, D. For $\,n>3\,$ there is only finitely additive 
rotationally invariant measure on the $\,n$-sphere on all Lebesgue 
measurable sets. Bull. AMS, {\bf{1}}, pp. 121-123, 1981.
\end{thebibliography}
\end{document}